\documentclass[12pt]{amsart}

\begin{document}

\title{Index transforms with Weber type kernels}

\author{S. Yakubovich}
\address{Department of Mathematics, Fac. Sciences of University of Porto,Rua do Campo Alegre,  687; 4169-007 Porto (Portugal)}
\email{ syakubov@fc.up.pt}

\keywords{Keywords  here} \subjclass[2000]{ 44A15, 33C10, 44A05}

\keywords{Index Transforms,  Bessel functions of the first and second kind,   Fourier transform,   Mellin transform,  Boundary value problem}

\maketitle

\markboth{\rm \centerline{ S.  Yakubovich}}{}
\markright{\rm \centerline{ Index transforms with Weber type kernels}}

\begin{abstract}   New index transforms with  Weber type kernels, consisting  of products of  Bessel  functions  of the first  and second kind   are  investigated.   Mapping properties and inversion formulas are established for these transforms  in Lebesgue  spaces.  The results are applied to solve a  boundary value problem on the wedge for a  fourth   order partial differential equation.  
\end{abstract}

\section{Introduction and preliminary results}

Let $\alpha \in \mathbb{R}$ and $ f(x),  \ g(\tau),\  x \in \mathbb{R}_+,\  \tau \in  \mathbb{R}$ be  complex-valued functions.  In this paper we will investigate mapping properties of the following index transforms   \cite{yak}
$$(F_\alpha f ) (\tau) =  \int_0^\infty  {\mathcal W}_\alpha(x,\tau) f(x)dx,  \eqno(1.1)$$
$$(G_\alpha  g) (x) =   \int_{-\infty}^\infty   {\mathcal W}_\alpha(x,\tau)\    g(\tau) d\tau, \eqno(1.2)$$
where ${\mathcal W}_\alpha(x,\tau)$ is the Weber type kernel defined by the formula

$${\mathcal W}_\alpha(x,\tau) =  { {\rm Im} \left[ J_{\alpha +i\tau} (\sqrt x)  Y_{\alpha - i\tau} (\sqrt x) \right]\over \sinh(\pi\tau)} .\eqno(1.3)$$
Here  $i$ is the imaginary unit, ${\rm Im}$ is the imaginary part of a complex number  and $J_\nu(z),\ Y_\nu(z)$ are  Bessel  functions  of the first and second kind, respectively,  \cite{erd}, Vol. II, which, in turn,  are fundamental solutions of the  differential equation
$$  z^2{d^2u\over dz^2}  + z{du\over dz} + (z^2- \nu^2)u = 0.\eqno(1.4)$$
Bessel functions in  (1.3) with fixed $\tau$ form the kernel of the well-known Weber-Orr integral transforms (cf. \cite{genweb}, \cite{mal},\ \cite{nasim} ) widely used, for instance, in the  heat process problems for  radial-symmetric
regions.   Concerning  index transforms of the Weber type, we mention  integral operators considered in \cite{vu}, where the variable of integration in the corresponding  kernel is involved in the arguments and indices of Bessel functions.  Appealing to relation ( 8.4.20.25) in \cite{prud}, Vol. III we find the important Mellin-Barnes integral representation of the kernel (1.3), namely,

$$ {\mathcal W}_\alpha(x,\tau) =  { \cosh(\pi\tau) \over 2\pi^{7/2} i } \int_{\gamma-i\infty}^{\gamma +i\infty}\Gamma(s+ i\tau)\Gamma(s-i\tau) $$

$$\times   \frac { \Gamma(s+\alpha)  \Gamma(1/2-s)\Gamma(1-s)}{ \Gamma (1+\alpha -s)  } x^{-s} ds, \ x >0,\eqno(1.5)$$
where $\Gamma(z)$ is Euler's gamma-function \cite{erd}, Vol. I  and $\gamma$ is chosen  from the interval $(\hbox{max}(-\alpha, 0),\ 1/2)$.    Integral (1.5) converges absolutely for any $\gamma$ from the open vertical strip $\hbox{max}(-\alpha, 0) <  \gamma < 1/2$ and uniformly with respect to $x \ge x_0 >0$ by virtue of the Stirling asymptotic formula for the gamma-function \cite{erd}, Vol. I, which yields for a fixed $\tau$

$$  \frac {\Gamma(s+ i\tau)\Gamma(s-i\tau)  \Gamma(s+\alpha)  \Gamma(1/2-s)\Gamma(1-s)}{ \Gamma (1+\alpha -s)  }$$

$$= O\left( e^{- 2\pi|s|} |s|^{ 2\gamma - 3/2}\right), \  |s| \to \infty.\eqno(1.6)$$

The Mellin-Barnes representation (1.5) can be used   to represent  ${\mathcal W}_\alpha(x,\tau)$  for all $x >0,\ \tau \in \mathbb{R}$ as the Fourier cosine transform \cite{tit}.   Precisely, we have

{\bf Lemma 1}. {\it Let $x >0,  \tau \in \mathbb{R},\  |\alpha| < 1/2, \ \alpha \neq 0$. Then $ {\mathcal W}_\alpha(x,\tau)$ has the following representation in terms of the Fourier cosine transform 
$${\mathcal W}_\alpha(x,\tau)  =   { 2\cosh(\pi\tau) \over \pi \sin(2\pi\alpha)}   \int_0^\infty  \left[  {\bf J}_{2\alpha} (2\sqrt x \cosh(u/2) ) \right.$$
$$\left.  -   J_{2\alpha} (2\sqrt x \cosh(u/2) ) \right]  \cos(\tau u) du,\eqno(1.7) $$
where ${\bf J}_\nu(z)$ is the Anger function  \cite{erd}, Vol. II.}

\begin{proof}     In fact,  appealing  to the reciprocal formulae via the Fourier cosine transform (cf. formula (1.104) in \cite{yak}) 
$$\int_0^\infty  \Gamma\left(s + i\tau\right)  \Gamma\left(s - i\tau \right)  \cos( \tau y) d\tau
= {\pi\over 2^{2s}}  {\Gamma(2s) \over \cosh^{2s}(y/2)},\ {\rm Re}\ s > 0,\eqno(1.8)$$
$$  \Gamma\left(s + i\tau \right)  \Gamma\left(s - i\tau\right)  
=   { \Gamma(2s)  \over 2^{2s-1}}  \int_0^\infty   {\cos(\tau y)  \over \cosh^{2s} (y/2)} \ dy,\eqno(1.9)$$ 
we replace  the gamma-product $ \Gamma\left(s + i\tau \right)  \Gamma\left(s - i\tau\right)$ in  (1.5)  by the right-hand side of  (1.9). Then, changing the order of integration via the absolute convergence, which can be justified using the Stirling asymptotic formula for the gamma-function, we employ the reflection and duplication formulae  for the gamma-function \cite{erd}, Vol. I to derive 

$$ {\mathcal W}_\alpha(x,\tau) =  { \cosh(\pi\tau) \over \pi^{2} i}    \int_0^\infty   \cos(\tau u)  \int_{\gamma-i\infty}^{\gamma +i\infty}   \frac { \Gamma(s+\alpha)  \left( x  \cosh^{2} (u/2)\right)^{-s}  }{ \sin(2\pi s) \ \Gamma (1+\alpha -s)  }  ds du.\eqno(1.10)$$
In the meantime, the inner integral with respect to $s$ in (1.10) is calculated in \cite{prud}, Vol. III, relation (8.4.26.5), and we have

$$   \int_{\gamma-i\infty}^{\gamma +i\infty}  \frac { \Gamma(s+\alpha)  \Gamma(s) \Gamma(s+1/2)\Gamma(1/2-s)\Gamma(1-s)}{ \Gamma (1+\alpha -s)  } x^{-s} ds$$

$$ = - {4\pi^3 i\over \sin(2\alpha\pi)} \left[ J_{2\alpha} (2\sqrt x) - {\bf J}_{2\alpha} (2\sqrt x) \right],$$
where  ${\bf J}_\nu(z)$ is the Anger function  \cite{erd}, Vol. II.  Hence, substituting this result in (1.10),  we get (1.7). 

\end{proof}

Further,  recalling  the Mellin-Barnes representation (1.5) of the kernel  ${\mathcal W}_\alpha(x,\tau)$, we will derive an ordinary differential equation whose particular solution is  ${\mathcal W}_\alpha(x,\tau)$. Precisely,   it is given by 

{\bf Lemma 2}.    {\it Let $\alpha < 1/2$.  For each $\tau \in \mathbb{R}$ the function   ${\mathcal W}_\alpha(x,\tau)$ satisfies the following fourth  order differential equation with variable  coefficients}
$$x^3 {d^4 {\mathcal W}_\alpha  \over dx^4} +   6 x^2 {d^3 {\mathcal W}_\alpha  \over dx^3} + x\left( 7+\tau^2-\alpha^2 +x\right) {d^2   {\mathcal W}_\alpha \over dx^2} $$

$$ +  \left( 1+\tau^2-\alpha^2 + {5\over 2}  x \right) {d {\mathcal W}_\alpha \over dx}  
+ \left( {1\over 2} - {(\alpha\tau)^2 \over x}  \right) {\mathcal W}_\alpha  = 0,\ x >0.\eqno(1.11)$$

\begin{proof}   Indeed,  the asymptotic behavior at infinity (1.6) of the integrand in (1.5) and the absolute and uniform convergence of the integral and its derivatives with respect to $x$ allow us to differentiate under the integral sign.     Hence, applying twice the differential operator $x{d\over dx} $ to both sides of (1.5) and employing the reduction formula for gamma-function \cite{erd}, Vol. I, we obtain  

$$ \left(x {d\over dx}\right)^2 {\mathcal W}_\alpha(x,\tau) =  { \cosh(\pi\tau) \over 2\pi^{7/2} i} \int_{\gamma-i\infty}^{\gamma +i\infty}s^2 \Gamma(s+ i\tau)\Gamma(s-i\tau) $$

$$\times   \frac { \Gamma(s+\alpha)  \Gamma(1/2-s)\Gamma(1-s)}{ \Gamma (1+\alpha -s)  } x^{-s} ds$$

$$=  { \cosh (\pi\tau)\over 2\pi^{7/2} i} \int_{\gamma-i\infty}^{\gamma +i\infty} \Gamma(s+ 1+ i\tau)\Gamma(s+1-i\tau)   \frac { \Gamma(s+\alpha)  \Gamma(1/2-s)\Gamma(1-s)}{ \Gamma (1+\alpha -s)  } x^{-s} ds$$

$$ - \tau^2 {\mathcal W}_\alpha(x,\tau) = { \cosh(\pi\tau)\over \pi^{7/2} } \int_{\gamma +1-i\infty}^{\gamma +1+i\infty} \Gamma(s+ i\tau)\Gamma(s-i\tau) $$

$$\times   \frac { \Gamma(s+\alpha-1)  \Gamma(3/2-s)\Gamma(2-s)}{ \Gamma (2+\alpha -s)  } x^{1-s} ds- \tau^2 {\mathcal W}_\alpha(x,\tau).$$
After multiplication by $x^{\alpha} $ and differentiation of both sides of the obtained equality, we employ again the reduction formula for the gamma-function to derive

$$ {d\over dx} \left( x^{\alpha} \left(x {d\over dx}\right)^2 {\mathcal W}_\alpha(x,\tau)  \right) = { \cosh(\pi\tau) \over 2\pi^{7/2} i} \int_{\gamma +1-i\infty}^{\gamma +1+i\infty} \Gamma(s+ i\tau)\Gamma(s-i\tau) $$

$$\times   \frac { \Gamma(s+\alpha-1)  \Gamma(3/2-s)\Gamma(2-s)}{ \Gamma (1+\alpha -s)  } x^{\alpha-s} ds- \tau^2{d\over dx} \left( x^\alpha  {\mathcal W}_\alpha(x,\tau)\right).\eqno(1.12)$$
In the meantime, the contour in the integral in the right-hand side of (1.12) can be moved to the left via Slater's theorem \cite{prud}, Vol. III. Hence, integrating over the vertical line ${\rm Re}\  s = \gamma$ in the complex plane, we take into account assumptions  $\hbox{max}(-\alpha, 0) < \gamma <  1/2,\  \alpha < 1/2$ to  guarantee the condition  $0< \gamma +\alpha < 1$.  Then the integral can be treated  as follows

$$\int_{\gamma -i\infty}^{\gamma +i\infty} \Gamma(s+ i\tau)\Gamma(s-i\tau) \   \frac { \Gamma(s+\alpha-1)  \Gamma(3/2-s)\Gamma(2-s)}{ \Gamma (1+\alpha -s)  } x^{\alpha-s} ds$$

$$=  \int_{\gamma -i\infty}^{\gamma +i\infty} \Gamma(s+ i\tau)\Gamma(s-i\tau)\    \frac {  (1/2-s) (1-s) \Gamma(s+\alpha) \Gamma(1/2-s)\Gamma(1-s)}{(s+\alpha-1) \Gamma (1+\alpha -s)  } x^{\alpha-s} ds$$

$$=   x^\alpha {d\over dx} \int_{\gamma -i\infty}^{\gamma +i\infty} \Gamma(s+ i\tau)\Gamma(s-i\tau)\    \frac {  (1/2-s)  \Gamma(s+\alpha) \Gamma(1/2-s)\Gamma(1-s)}{(s+\alpha-1) \Gamma (1+\alpha -s)  } x^{1-s} ds$$

$$ = x^\alpha {d\over dx} \  x^{3/2} {d\over dx} \int_{\gamma -i\infty}^{\gamma +i\infty} \Gamma(s+ i\tau)\Gamma(s-i\tau)\    \frac { \Gamma(s+\alpha) \Gamma(1/2-s)\Gamma(1-s)}{(s+\alpha-1) \Gamma (1+\alpha -s)  } x^{ 1/2 -s} ds$$

$$ = - x^\alpha {d\over dx} \  x^{3/2} {d\over dx}  x^{\alpha- 1/2} \int_0^x \int_{\gamma -i\infty}^{\gamma +i\infty} \Gamma(s+ i\tau)\Gamma(s-i\tau) $$

$$\times \frac { \Gamma(s+\alpha) \Gamma(1/2-s)\Gamma(1-s)}{ \Gamma (1+\alpha -s)  } y^{ -s-\alpha} ds dy.$$
Substituting this result into (1.12) and using (1.5),  we find

$$ {d\over dx} \left( x^{\alpha} \left(x {d\over dx}\right)^2 {\mathcal W}_\alpha(x,\tau)  \right) +  \tau^2{d\over dx} \left( x^\alpha  {\mathcal W}_\alpha(x,\tau)\right) $$

$$=  - x^\alpha {d\over dx} \  x^{3/2} {d\over dx}  x^{\alpha- 1/2} \int_0^x  y^{-\alpha}  {\mathcal W}_\alpha(y,\tau) dy.\eqno(1.13)$$
Finally, fulfilling the differentiation in (1.13), we end up with equation (1.11), completing the proof of Lemma 2. 

\end{proof} 

 In the sequel we will employ the Mellin  transform technique developed in \cite{yal} in order to  investigate the mapping properties of the index transforms  (1.1), (1.2).   Precisely, the Mellin transform is defined, for instance, in  $L_{\nu, p}(\mathbb{R}_+),\ 1 \le  p \le 2$ (see details in \cite{tit}) by the integral  
$$f^*(s)= \int_0^\infty f(x) x^{s-1} dx,\eqno(1.14)$$
 being convergent  in mean with respect to the norm in $L_q(\nu- i\infty, \nu + i\infty),\ \nu \in \mathbb{R}, \   q=p/(p-1)$.   Moreover, the  Parseval equality holds for $f \in L_{\nu, p}(\mathbb{R}_+),\  g \in L_{1-\nu, q}(\mathbb{R}_+)$
$$\int_0^\infty f(x) g(x) dx= {1\over 2\pi i} \int_{\nu- i\infty}^{\nu+i\infty} f^*(s) g^*(1-s) ds.\eqno(1.15)$$
The inverse Mellin transform is given accordingly
 $$f(x)= {1\over 2\pi i}  \int_{\nu- i\infty}^{\nu+i\infty} f^*(s)  x^{-s} ds,\eqno(1.16)$$
where the integral converges in mean with respect to the norm  in   $L_{\nu, p}(\mathbb{R}_+)$
$$||f||_{\nu,p} = \left( \int_0^\infty  |f(x)|^p x^{\nu p-1} dx\right)^{1/p}.\eqno(1.17)$$
In particular, letting $\nu= 1/p$ we get the usual space $L_p(\mathbb{R}_+; \  dx)$.

\section {Boundedness  and inversion properties of the index transform (1.1)}

Let $C_b (\mathbb{R};  [ \cosh(\pi\tau)]^{-1})$ be the space of functions $f(\tau)$ such that $f(\tau)/  \cosh (\pi\tau) $ is bounded continuous on $\mathbb{R}$.  

{\bf Theorem 2.}   {\it Let $\hbox{max} (-\alpha, 0) < \gamma <  1/2$. The index transform  $(1.1)$  is well-defined as a  bounded operator $F_\alpha: L_{1-\gamma,1} \left(\mathbb{R}_+\right) \to C_b (\mathbb{R};  [ \cosh (\pi\tau)]^{-1})$ and the following norm inequality takes place
$$||F_\alpha f||_{C_b} \equiv \sup_{\tau \in \mathbb{R}}  {(\left| F_\alpha f)(\tau)\right| \over  \cosh (\pi\tau)} \le C_\gamma  ||f||_{L_{1-\gamma,1}\left(\mathbb{R}_+\right)},\eqno(2.1)$$
where 

$$C_\gamma =  {2^{2\gamma -1}   \over \pi^{2} }  B(\gamma, \gamma)   \int_{\gamma-i\infty}^{\gamma +i\infty}   \frac { \left|\Gamma(s+\alpha)\right|  }{ \left|\sin(2\pi s) \ \Gamma (1+\alpha -s) \right| }  |ds |\eqno(2.2)$$
and $B(a,b)$ is the beta-function \cite{erd}, Vol. I.   Moreover, it admits the integral representation 

$$(F_\alpha f) (\tau) = { 1\over  \pi^2}\int_0^\infty K_{i\tau} \left(\sqrt x\right) \left[ I_{i\tau} \left(\sqrt x\right) + I_{-i\tau} \left(\sqrt x\right)\right] \varphi_\alpha(x) dx, \ \tau \in \mathbb{R},\eqno(2.3)$$
where

$$ \varphi_\alpha(x)  = {1\over 2\pi i}   \int_{1-\gamma - i\infty}^{1-\gamma +i\infty}  \frac{ \Gamma(1-s+\alpha) \Gamma^2(s)}{\Gamma(s+\alpha)}  f^*(s) x^{-s} ds,\  x >0\eqno(2.4)$$
and integrals $(2.3)$, $(2.4)$ converge absolutely.}

\begin{proof}  In fact, from (1.10) we have the estimate 

$$ \left|{\mathcal W}_\alpha(x,\tau)\right| \le  { 2 x^{-\gamma}  \over \pi^{2} } \  \cosh (\pi\tau)  \int_0^\infty   {du\over \cosh^{2\gamma} (u) }   \int_{\gamma-i\infty}^{\gamma +i\infty}   \frac { \left|\Gamma(s+\alpha)\right|  }{ \left|\sin(2\pi s) \ \Gamma (1+\alpha -s) \right| }  |ds |$$

$$= C_\gamma x^{-\gamma}  \cosh (\pi\tau),\eqno(2.5)$$
where the constant $C_\gamma$ is defined by (2.2).   Hence,

$$ ||F_\alpha f ||_{C_b} = \sup_{\tau \in \mathbb{R}} {\left| (F_\alpha f)(\tau)\right|\over  \cosh (\pi\tau)}  \le
C_\gamma  \int_0^\infty  x^{-\gamma} |f(x)| dx,$$
which yields (2.1) via (1.17).  Next, recalling Parseval equality (1.15), the Mellin-Barnes representation (1.5) of the kernel  ${\mathcal W}_\alpha(x,\tau)$ and formula (1.16) of the Mellin transform, we have 

$$(F_\alpha f) (\tau) =  {\cosh (\pi\tau) \over 2\pi^{7/2} i} \int_{\gamma-i\infty}^{\gamma +i\infty}\Gamma(s+ i\tau)\Gamma(s-i\tau) $$

$$\times   \frac { \Gamma(s+\alpha)  \Gamma(1/2-s)\Gamma(1-s)}{ \Gamma (1+\alpha -s)  } f^*(1-s) ds.\eqno(2.6)$$
Meanwhile, appealing to relation (8.4.23.23) in \cite{prud}, Vol. III,  it implies the representation 

$${\sqrt \pi \over \cosh(\pi\tau)} K_{i\tau} \left(\sqrt x\right) \left[ I_{i\tau} \left(\sqrt x\right) + I_{-i\tau} \left(\sqrt x\right)\right]=
{1\over 2\pi i} \int_{\gamma-i\infty}^{\gamma +i\infty}\Gamma(s+ i\tau)\Gamma(s-i\tau) $$

$$\times   \frac {  \Gamma(1/2-s)}{\Gamma(1-s)} x^{-s} ds, \ x> 0,\ 0< \gamma < {1\over 2}.$$  
Hence employing again Parseval equality (1.15), it can be rewritten in the form (2.3), where integrals (2.3), (2.4) converge absolutely owing to Stirling asymptotic formula for the gamma-function,  the asymptotic behavior of the kernel in (2.3) (see \cite{erd}, Vol. II)

$$K_{i\tau} \left(\sqrt x\right) \left[ I_{i\tau} \left(\sqrt x\right) + I_{-i\tau} \left(\sqrt x\right)\right] =O\left( \log x\right),\ x \to 0,$$

$$K_{i\tau} \left(\sqrt x\right) \left[ I_{i\tau} \left(\sqrt x\right) + I_{-i\tau} \left(\sqrt x\right)\right] =O\left( {1\over \sqrt x}\right),\ x \to \infty,$$
and the condition  $\hbox{max} (-\alpha, 0) < \gamma <  1/2$.   Theorem 2 is proved. 

\end{proof}

Now we will derive a composition relation of the index transform (1.1) with the Kontorovich-Lebedev transform \cite{yak}. In fact,  taking relation (8.4.23.5) in \cite{prud}, Vol.III,  it reads 

$${\sqrt\pi\over \cosh(\pi\tau) } e^{x/2} K_{i\tau}\left({x\over 2} \right) = {1\over 2\pi i} \int_{\gamma-i\infty}^{\gamma +i\infty}\Gamma(s+ i\tau)\Gamma(s-i\tau) $$

$$\times  \Gamma(1/2-s) x^{-s} ds, \ x> 0,\ 0< \gamma < {1\over 2}.\eqno(2.7)$$  
Therefore, recalling (2.6), we write via the same arguments the representation of $(F_\alpha f) (\tau)$ in terms of the Kontorovich-Lebedev transform 

$$(F_\alpha f) (\tau) = {2\over \pi^2} \int_0^\infty K_{i\tau} (x) e^x \psi_\alpha(2x) dx,\eqno(2.8)$$
where 

$$\psi_\alpha(x)  = {1\over 2\pi i}   \int_{1-\gamma - i\infty}^{1-\gamma +i\infty}  \frac{ \Gamma(1-s+\alpha) \Gamma(s)}{\Gamma(s+\alpha)}  f^*(s) x^{-s} ds,\  x >0.\eqno(2.9)$$
If $e^x \psi_\alpha(2x)  \in L_{1,2}\left(\mathbb{R}_+\right)$, then integral (2.8) converges in the mean square sense with respect to the norm in $L_2\left(\mathbb{R}_+; \tau\sinh(\pi\tau) d\tau\right)$ (cf. \cite{yak}),  $(F_\alpha f) (\tau) \in  L_2\left(\mathbb{R}_+; \tau\sinh(\pi\tau) d\tau\right)$ and the Parseval equality holds

$$\int_0^\infty  e^{x} \left|\psi_\alpha(x)\right|^2  x dx =  4 \int_0^\infty \tau \sinh(\pi\tau)  \left| (F_\alpha f) (\tau)  \right|^2 d\tau.\eqno(2.10)$$
 Moreover, the inversion formula for the Kontorovich-Lebedev transform implies 

 $$ x e^x \psi_\alpha(2x) =  \int_0^\infty \tau \sinh(\pi\tau)  K_{i\tau} (x) (F_\alpha f) (\tau) d\tau,\eqno(2.11)$$
 where the integral converges in the mean square sense with  respect to the norm in $L_{1,2}\left(\mathbb{R}_+\right)$. 
Identity (2.11) will be a starting point to prove the inversion formula for the index transform (1.1).  In fact, taking into account (2.9), we get from (2.11)

$$   {1\over 2\pi i}   \int_{1-\gamma - i\infty}^{1-\gamma +i\infty}  \frac{ \Gamma(1-s+\alpha) \Gamma(s)}{\Gamma(s+\alpha)}  f^*(s) x^{1-s} 2^{-s} ds$$

$$ =  \int_0^\infty \tau \sinh(\pi\tau)  e^{-x} K_{i\tau} (x) (F_\alpha f) (\tau) d\tau.\eqno(2.12)$$
Then relation (8.4.23.3) in \cite{prud}, Vol. III gives the integral 

$$ {1\over \sqrt\pi } \  e^{- x/2} K_{i\tau}\left({x\over 2} \right) = {1\over 2\pi i} \int_{\nu-i\infty}^{\nu +i\infty}{\Gamma(s+ i\tau)\Gamma(s-i\tau) \over  \Gamma(1/2+ s)}  x^{-s} ds, \ x> 0, \ \nu > 0.\eqno(2.13)$$  
Hence, making simple substitutions and changes of variables in (2.12), we take  the Mellin transform of both sides over some line ${\rm Re} s=\nu >0$ in the analyticity strip $a< {\rm Re} s < b,\ b >0$ of $f^*(s)$ to obtain 

$$    \frac{ \Gamma(\alpha-s) \Gamma(s+1)}{\Gamma(1+s+\alpha)}  f^*(1+s)  =   2 \int_0^\infty x^{s-1} \int_0^\infty \tau \sinh(\pi\tau) (F_\alpha f) (\tau)$$

$$\times  e^{-x/2} K_{i\tau} \left({x\over 2}\right)  d\tau dx,\ \alpha > 0,\eqno(2.14)$$
where $ \hbox{max} (a-1, 0) < \nu < \hbox{min} (\alpha, b-1)$. This strip can be easily established via the absolute convergence of the corresponding integrals if we assume that $f$ be locally integrable on $\mathbb{R}_+$, i.e. $f(x)\in L_{loc} (\mathbb{R}_+)$ and behaves as  $f(x)= O\left(x^{-a}\right),\ x \to 0,\ f(x)= O\left(x^{-b}\right),\ x \to \infty.$   To change the order of integration in the right-hand side of (2.14), we use the following estimate, basing on the Cauchy- Schwarz inequality and relation (2.16.32.2) in \cite{prud}, Vol. II

$$ \int_0^\infty x^{\nu -1} \int_0^\infty \tau \sinh(\pi\tau) \left|(F_\alpha f) (\tau)\right|  e^{-x/2} \left|K_{i\tau} \left({x\over 2}\right) \right|  d\tau dx $$

$$ \le \left(\int_0^\infty \tau (1+\tau)^2  \sinh(\pi\tau) \left|(F_\alpha f) (\tau)\right|^2 \right)^{1/2} $$

$$\times  \left(\int_0^\infty {\tau \sinh(\pi\tau)\over (1+\tau)^2} \left(\int_0^\infty x^{\nu -1}  e^{-x/2} \left|K_{i\tau} \left({x\over 2}\right) \right| dx\right)^2 d\tau\right)^{1/2} $$

$$ \le 2^{\nu/2} \left[\Gamma(\nu )\right]^{1/2}  \left(\int_0^\infty \tau (1+\tau)^2  \sinh(\pi\tau) \left|(F_\alpha f) (\tau)\right|^2 \right)^{1/2} $$

$$\times  \left(\int_0^\infty {\tau \sinh(\pi\tau)\over (1+\tau)^2} \int_0^\infty x^{\nu -1}  K^2_{i\tau} \left(x\right) dx d\tau\right)^{1/2} $$

$$=  2^{\nu -3/2}   \Gamma\left({\nu\over 2}\right) \left(\int_0^\infty \tau (1+\tau)^2  \sinh(\pi\tau) \left|(F_\alpha f) (\tau)\right|^2 \right)^{1/2} $$

$$\times  \left(\int_0^\infty {\tau \sinh(\pi\tau)\over (1+\tau)^2} \left|\Gamma\left(  {\nu\over 2}  +i\tau\right)\right|^2 d\tau\right)^{1/2} < \infty,  \  0 < \nu <  1 $$
under the additional condition 

$$(F_\alpha f) (\tau) \in L_2\left(\mathbb{R}_+; \tau (1+\tau)^2  \sinh(\pi\tau)d\tau \right) \subset 
L_2\left(\mathbb{R}_+; \tau \sinh(\pi\tau)d\tau \right).$$
Hence, returning to (2.14) and calculating the integral with respect to $x$ via relation (8.4.23.3),  we multiply both sides by  the product $\Gamma(1/2-s)\Gamma(1-s)$ and with the aid of the duplication and reflection formulae for the gamma-function write it in the form 

$$  { s f^*(1+s) \over \sin(2\pi s) (\alpha^2-s^2)}   =  \frac{ \Gamma(\alpha+s) \Gamma(1/2-s)\Gamma(1-s)}{\pi\sqrt \pi \Gamma(1-s+\alpha)} \int_0^\infty \tau \sinh(\pi\tau) (F_\alpha f) (\tau)$$

$$\times \Gamma(s+ i\tau) \Gamma(s- i\tau)   d\tau,\  0< {\rm Re} s < {1\over 2}.\eqno(2.15)$$
Now  we take $\nu$ from the interval  $( 0, 1/2)$ and integrate both sides over $(\nu-i\infty, \nu+i\infty)$.  Interchanging the order of integration in the obtained equality, we appeal to (1.5) to derive

$$   {1\over 2\pi i}  \int_{\nu-i\infty}^{\nu +i\infty} { s f^*(1+s) \over \sin(2\pi s) (\alpha^2-s^2)} x^{-s} ds $$

$$  =  \pi  \int_0^\infty \tau \tanh(\pi\tau) {\mathcal W}_\alpha(x,\tau) (F_\alpha f) (\tau) d\tau.\eqno(2.16)$$
To justify the interchange of the order of integration, we employ Fubini's theorem, the Stirling asymptotic formula for the gamma-function and the following estimate

$$\int_0^\infty \tau \sinh(\pi\tau) \left| (F_\alpha f) (\tau) \Gamma(s+ i\tau) \Gamma(s- i\tau) \right|   d\tau$$

$$\le \left(\int_0^\infty \tau^2  \left|(F_\alpha f) (\tau)\right|^2 d\tau\right)^{1/2} \left(\int_0^\infty \sinh^2 (\pi\tau) \left|\Gamma(s+ i\tau) \Gamma(s- i\tau) \right|^2   d\tau\right)^{1/2}$$
$$\le {1\over \sqrt\pi} \left(\int_0^\infty \tau (1+\tau)^2  \sinh(\pi\tau) \left|(F_\alpha f) (\tau)\right|^2 d\tau\right)^{1/2} $$

$$\times \left(\int_0^\infty \sinh^2 (\pi\tau) \left|\Gamma(s+ i\tau) \Gamma(s- i\tau) \right|^2   d\tau\right)^{1/2}.\eqno(2.17)$$
The right-hand side of the second  inequality in (2.17) is finite if we show the upper bound for the latter integral.   However, it can be calculated explicitly. Indeed,  employing relation (2.5.46.15) in \cite{prud}, Vol. I,  we find the value of the following Fourier sine transform  

$$\int_0^\infty {\sin(t\tau)\over \sinh^{2s} (t/2)} dt = {\sinh(\pi\tau)\over \sqrt \pi}  \frac{ \Gamma(s+ i\tau) \Gamma(s- i\tau) \Gamma(1-s)}{\Gamma(1/2 + s)},\eqno(2.18)$$
where $ \tau >0$ and the integral converges absolutely under condition $ 0 < {\rm Re} s < 1/2$.  When $ 0 < {\rm Re}  s < 1/4$ the hyperbolic function $\sinh^{-2s} (t/2) \in L_2\left(\mathbb{R}_+\right)$, and the reciprocal formula for the Fourier sine transform \cite{tit} says that the right-hand side of (2.17) belongs to $L_2\left(\mathbb{R}_+\right)$ and 

$$\int_0^\infty  \sinh(\pi\tau) \Gamma(s+ i\tau) \Gamma(s- i\tau)  \sin(t \tau ) d\tau = {\pi\sqrt \pi \ \Gamma(1/2 + s) \over 2 \Gamma(1-s) \sinh^{2s} (t/2)},\eqno(2.19)$$ 
where the integral converges in the mean square sense when $ 0 < {\rm Re s} < 1/4$ and relatively under condition $ 0 < {\rm Re s} < 1/2$. Besides, the Parseval  identity implies the equality

$$\int_0^\infty  \sinh^2(\pi\tau) \left| \Gamma(s+ i\tau) \Gamma(s- i\tau) \right|^2  d\tau$$

$$ = \pi^2 \left| {\Gamma(1/2 + s) \over  \Gamma(1-s)}\right|^2 \int_0^\infty {dt \over  \sinh^{4{\rm Re} s} (t)},\eqno(2.20)$$
where $ 0 < {\rm Re} s < 1/4$.  But the integral in the right-hand side of (2.19) is calculated via relation (2.5.46.15) (slightly corrected) in \cite{prud}, Vol. I, and we obtain a possibly new formula

$$\int_0^\infty  \sinh^2(\pi\tau) \left| \Gamma(s+ i\tau) \Gamma(s- i\tau) \right|^2  d\tau $$

$$= {\pi \sqrt\pi \over 2} \Gamma\left(1/2- 2{\rm Re} s\right)\Gamma\left(2{\rm Re} s \right)  \left| {\Gamma(1/2 + s) \over  \Gamma(1-s)}\right|^2,\ 0 < {\rm Re} s < {1\over 4}.\eqno(2.21)$$
Returning to the equality (2.16), we restrict $\nu$ for the interval $\hbox{max} (a-1, 0) < \nu < \hbox{min} (1/4, \alpha, b-1)$ and apply to its both sides the differential operator $x {d\over dx}$. Then differentiating under the integral sign in the left-hand side  owing to the boundedness of $f^*(1+s)$ in the analyticity strip, we find

$$   {1\over 2\pi i}  \int_{\nu -i\infty}^{\nu +i\infty} {  f^*(1+s) \over \sin(2\pi s) } x^{-s} ds + {\alpha^2 \over 2\pi i}  \int_{\nu-i\infty}^{\nu +i\infty} {  f^*(1+s) \over \sin(2\pi s) (s^2-\alpha^2)} x^{-s} ds $$

$$  =  \pi  x {d\over dx} \int_0^\infty \tau \tanh(\pi\tau) {\mathcal W}_\alpha(x,\tau) (F_\alpha f) (\tau) d\tau.\eqno(2.22)$$
Acting again by the differential operator $x {d\over dx}$ on (2.22) and using (2.16), we derive the equality 

$$    {1\over 2\pi i}  \int_{\nu-i\infty}^{\nu +i\infty} {   s f^*(1+s) \over \sin(2\pi s) } x^{-s} ds  $$

$$  =  \pi  \left[\alpha^2- \left(x {d\over dx}\right)^2 \right] \int_0^\infty \tau \tanh(\pi\tau) {\mathcal W}_\alpha(x,\tau) (F_\alpha f) (\tau) d\tau.\eqno(2.23)$$
Continuing this process, we use the infinite product for the sine function, namely,

$${\sin(2\pi s)\over 2\pi s}= \prod_{n=1}^\infty \left(1- {4s^2\over n^2} \right),$$
and the condition $f^*(1+s) \in L_1(\nu-i\infty,  \nu +i\infty)$ to derive the following operator equality 

$$  xf(x)=  -i  \prod_{n=1}^\infty \left(1- {4\over n^2} \left(x {d\over dx}\right)^2  \right)  \left[   \int_{\nu-i\infty}^{\nu +i\infty} {   s f^*(1+s) \over \sin(2\pi s) } x^{-s} ds \right].$$
Hence, comparing with (2.23), we establish the inversion formula for the index transform (1.1)

$$x f(x)= 2\pi^2 \lim_{N\to \infty }  \int_0^\infty \tau \tanh(\pi\tau) $$

$$\times  \prod_{n=1}^N  \left(1- {4\over n^2} \left(x {d\over dx}\right)^2  \right)  \left(\alpha^2- \left(x {d\over dx}\right)^2 \right)  \left[ {\mathcal W}_\alpha(x,\tau)\right]  (F_\alpha f) (\tau) d\tau.\eqno(2.24)$$ 
We note that the differentiation under the integral sign in the right-hand side of (2.24) is allowed under the imposed conditions on $(F_\alpha f)(\tau)$.   

We summarize these results in the following inversion theorem for the transform (1.1).

{\bf Theorem 3}. {\it Let $\alpha > 0,   f$ be locally integrable on $\mathbb{R}_+$ and behave as  $f(x)= O\left(x^{-a}\right),\ x \to 0,\ f(x)= O\left(x^{-b}\right),\ x \to \infty,\ a< b, b > 1, a < 1+\alpha$.  If  $f^*(1+s) \in L_1(\nu-i\infty,  \nu +i\infty)$ for some $\nu$, satisfying the condition $\hbox{max} (a-1, 0) < \nu < \hbox{min} (1/4, \alpha, b-1)$ and $(F_\alpha f) (\tau) \in L_2\left(\mathbb{R}_+; \tau (1+\tau)^2  \sinh(\pi\tau)d\tau \right) $, then for all $x >0$ the inversion formula $(2.24)$ holds for the index transform $(1.1)$. }

In the meantime, relation (8.4.20.21) in \cite{prud}, Vol. III gives the Mellin-Barnes representation for the following Weber type kernel

$${\rm Re} \left[ J_{\alpha +i\tau} (\sqrt x)  Y_{\alpha - i\tau} (\sqrt x) \right] $$

$$= -   {1\over 2\pi \sqrt\pi i}   \int_{\nu - i\infty}^{\nu +i\infty}  \frac{ \Gamma(s+ i\tau) \Gamma(s- i\tau) \Gamma(\alpha+s)}{\Gamma(s) \Gamma(1/2 + s)\Gamma(1+\alpha-s)} x^{-s} ds,\eqno(2.25)$$
where the integral converges  for each $x >0$ and $\tau \in \mathbb{R}$ under the condition $ -\alpha < \nu < 3/4$.  Moreover, using the straightforward integration, we find the formula 

$$ -   {1\over 2\pi \sqrt\pi i}   \int_{\nu - i\infty}^{\nu +i\infty}  \frac{ \Gamma(s+ i\tau) \Gamma(s- i\tau) \Gamma(\alpha+s)}{\Gamma(1+s) \Gamma(1/2 + s)\Gamma(1+\alpha-s)} x^{-s} ds$$

$$ = \int_x^\infty  {\rm Re} \left[ J_{\alpha +i\tau} (\sqrt y)  Y_{\alpha - i\tau} (\sqrt y) \right] {dy\over y}.\eqno(2.26)$$
Hence we have

{\bf Remark 1}. {\it If  one pass to the limit  under the integral sign in $(2.24)$, this formula can be written as follows}

$$ f(x)= - {2\pi\over x}    \left(\alpha^2- \left(x {d\over dx}\right)^2 \right)  \int_0^\infty \int_x^\infty \tau  {\rm Re} \left[ J_{\alpha +i\tau} (\sqrt y)  Y_{\alpha - i\tau} (\sqrt y) \right] $$

$$\times  (F_\alpha f) (\tau)  {dy d\tau \over y}.\eqno(2.27)$$ 

Letting in (1.5) $\alpha=0$ and appealing to relation (8.4.20.35) in \cite{prud}, Vol. III,  operator (1.1) turns to be an  adjoint operator of the index transform with  Nicholson's  function as the kernel (cf. \cite{yaknic}) 

$$(F_0 f)(\tau)  = {1\over 2} \int_0^\infty \left[ J^2_{i\tau} \left(\sqrt x\right) +  Y^2_{i\tau} \left(\sqrt x\right)\right] f(x) dx.\eqno(2.28)$$

{\bf Corollary 1}. {\it Let conditions of Theorem 3 hold and $\alpha = 0$. Then for all $x >0$ the inversion formula holds 

$$x f(x)= - \pi^2 \lim_{N\to \infty }  \int_0^\infty \tau \tanh(\pi\tau)   \prod_{n=1}^N  \left(1- {4\over n^2} \left(x {d\over dx}\right)^2  \right)  \left(x {d\over dx}\right)^2  \left[ J^2_{i\tau} \left(\sqrt x\right)\right.$$

$$\left.  +  Y^2_{i\tau} \left(\sqrt x\right)\right]  (F_0 f) (\tau) d\tau.\eqno(2.29)$$ 
If, besides, $(F_0 f) (\tau) \in L_1\left(\mathbb{R}_+; \tau \cosh (\pi\tau) d\tau \right) $, it takes the form

$$ f(x) = - \pi  {d\over dx}  \int_0^\infty    \tau  {\rm Im} \left[  J^2_{i\tau} (\sqrt x )\right]  F(\tau)  d\tau,\eqno(2.30)$$
where the integral converges absolutely.} 

\begin{proof} Indeed, formula (2.29) is the corresponding equality (2.24) with $\alpha =0$.  On the other hand, combining with (2.13), relation (8.4.19.19) in \cite{prud}, Vol. III and Parseval equality (1.15), we deduce the integral representation

$$ {1\over  \pi}   \sinh(\pi\tau) \  e^{- x/2} K_{i\tau}\left({x\over 2} \right) = {1\over x} \int_0^\infty e^{- t/x} \ {\rm Im} \left[  J^2_{i\tau} (\sqrt t )\right] dt,\ x >0. \eqno(2.31)$$
Meanwhile, one can put $\alpha=0$ in  (2.14) after integration by parts in its right-hand side.  Hence

$$    \Gamma(1-s) f^*(1+s)  = 2  \int_0^\infty x^{s} {d\over dx} \int_0^\infty \tau \sinh(\pi\tau) (F_0 f) (\tau)$$

$$\times  e^{-x/2} K_{i\tau} \left({x\over 2}\right)  d\tau dx.\eqno(2.32)$$
Taking the inverse Mellin transform of both sides in (2.32), using (1.15) and differentiation under the integral sign, we have

$${d\over dx} \int_0^\infty e^{- t/x}  f(t) dt =  2  {d\over dx} \int_0^\infty \tau \sinh(\pi\tau) e^{-x/2} K_{i\tau} \left({x\over 2}\right) (F_0 f) (\tau) d\tau.\eqno(2.33)$$
Employing (2.30) and integrating by $x$ in (2.33), we find

$$ \int_0^\infty e^{- t/x} f(t) dt =   {2 \pi \over x} \int_0^\infty \int_0^\infty e^{- t/x} \ \tau {\rm Im} \left[  J^2_{i\tau} (\sqrt t )\right]  (F_0 f) (\tau) dt d\tau.\eqno(2.34)$$
After integration by parts in the left-hand side of (2.34) it becomes 

$$ \int_0^\infty e^{- t/x} \int_0^t f(u) du dt =   2 \pi  \int_0^\infty \int_0^\infty e^{- t/x} \ \tau {\rm Im} \left[  J^2_{i\tau} (\sqrt t )\right]  (F_0 f) (\tau) dt d\tau.\eqno(2.35)$$
But as  it follows from the Poisson integral for Bessel functions 

$$ \left| {\rm Im} \left[  J^2_{i\tau} (\sqrt t )\right]\right| \le \cosh(\pi\tau),\ t \ge 0.$$
Therefore if $(F_0 f) (\tau) \in L_1\left(\mathbb{R}_+; \tau \cosh (\pi\tau)d\tau \right) $, the iterated integral in the right-hand side of (2.35) converges absolutely. Hence the change of the order of integration is permitted. Finally, the injectivity of the Laplace transform drives to formula (2.29) after the corresponding differentiation.  Hence Corollary 1 is proved.  

\end{proof}

\section{Index transform (1.2)} 

In this section we will study the  boundedness  properties and  prove an inversion formula for  the index transform  (1.2).   We begin with 

{\bf Theorem 4.}  {\it  Let $\hbox{max} (-\alpha, 0) < \gamma <  1/4$ and  $ g(\tau)  \in L_2(\mathbb{R};  d\tau)$. Then $x^{\gamma} (G_\alpha g)(x)$ is bounded continuous on $\mathbb{R}_+$ and it holds  
$$\sup_{x >0}   \left| x^{\gamma}  (G_\alpha g)(x) \right|   \le C ||g ||_{L_2(\mathbb{R};  d\tau)},\eqno(3.1)$$
where 

$$C=  { ||g ||_{L_2(\mathbb{R};  d\tau)}\over 2 \pi^{7/4}} \left[  \Gamma\left(1/2- 2{\rm Re} s\right)\Gamma\left(2{\rm Re} s \right) \right]^{1/2} \int_{\gamma-i\infty}^{\gamma +i\infty}   \left | \frac { \Gamma(s+\alpha)}{ \sin(\pi s) \Gamma (1+\alpha -s)  } ds\right|.$$
 Besides, if 
 $$(G_\alpha g) (x) \in L_{\gamma,1}(\mathbb{R}_+)$$
 and its Mellin transform $|s|^{1/2-\gamma} e^{\pi|s|/2} (G_\alpha g)^*(s) \in L_1( \gamma - i\infty,\  \gamma +i\infty)$, then for all $y >0$ }
$${1\over 2\pi i}   \int_{\gamma -i\infty}^{\gamma  +i\infty} \frac{\Gamma(1+\alpha-s)}{\Gamma(\alpha+s)\Gamma(1-s)} (G_\alpha g)^*(s)  y^{ -s} ds$$
$$= {1\over \pi^2}   \int_{-\infty}^\infty  e^{y/2} \   K_{i\tau} \left({y\over 2} \right)  g(\tau)  d\tau. \eqno(3.2)$$

\begin{proof}   In fact, recalling (1.5) and employing the Cauchy-Schwarz and generalized Minkowski inequalities,  we find the estimate

$$  \left|  (G_\alpha g)(x) \right|  \le {1\over 2\pi^{7/2}} \int_{-\infty}^\infty  \cosh(\pi\tau) |g(\tau)|  \left|\int_{\gamma-i\infty}^{\gamma +i\infty} \Gamma(s+ i\tau)\Gamma(s-i\tau) \right.$$

$$\times  \left. \frac { \Gamma(s+\alpha)  \Gamma(1/2-s)\Gamma(1-s)}{ \Gamma (1+\alpha -s)  } x^{-s} ds\right| d\tau 
\le { ||g ||_{L_2(\mathbb{R};  d\tau)}\over 2\pi^{7/2}} $$

$$\times \left( \int_{-\infty}^\infty  \cosh^2(\pi\tau) \left| \int_{\gamma-i\infty}^{\gamma +i\infty}\Gamma(s+ i\tau)\Gamma(s-i\tau) \frac { \Gamma(s+\alpha)  \Gamma(1/2-s)\Gamma(1-s)}{ \Gamma (1+\alpha -s)  } x^{-\gamma} ds\right|^2 d\tau\right)^{1/2} $$

$$\le { x^{-\gamma}  ||g ||_{L_2(\mathbb{R};  d\tau)}\over 2\pi^{7/2}} \int_{\gamma-i\infty}^{\gamma +i\infty}  \left( \int_{-\infty}^\infty  \cosh^2(\pi\tau) \left|\Gamma(s+ i\tau)\Gamma(s-i\tau) \right|^2 d\tau\right)^{1/2} $$

$$\times \left | \frac { \Gamma(s+\alpha)  \Gamma(1/2-s)\Gamma(1-s)}{ \Gamma (1+\alpha -s)  } ds\right|.\eqno(3.3)$$
Our goal now is to prove an analog of formula (2.21) in order to estimate the latter iterated integral (3.3).   First we observe via the Stirling asymptotic formula for the gamma-function that the integral 

$$ \int_{-\infty}^\infty  \cosh^2(\pi\tau) \left|\Gamma(s+ i\tau)\Gamma(s-i\tau) \right|^2 d\tau$$
converges for any $s$ form the strip $0 < {\rm Re} s < 1/4.$ Then appealing to the reciprocal formula for the Fourier cosine transform, which corresponds  to the slightly corrected relation (2.5.46.15) in \cite{prud}, Vol. I, we arrive at the following analog of formula (2.19)

$$  \int_{0}^\infty  \cosh(\pi\tau) \Gamma(s+ i\tau)\Gamma(s-i\tau) \cos(t\tau) d\tau = {\pi\sqrt \pi \ \Gamma( s) \over 2 \Gamma(1/2-s) \sinh^{2s} (t/2)}.\eqno(3.4)$$ 
This integral converges in the mean square sense when $0< {\rm Re} s < 1/4$ and relatively for $0< {\rm Re} s < 1/2$.  Hence via Parseval's equality we find

$$\int_{-\infty}^\infty  \cosh^2(\pi\tau) \left|\Gamma(s+ i\tau)\Gamma(s-i\tau) \right|^2 d\tau $$

$$=  \pi \sqrt\pi \Gamma\left(1/2- 2{\rm Re} s\right)\Gamma\left(2{\rm Re} s \right)  \left|{ \Gamma( s) \over  \Gamma(1/2-s)} \right|^2 ,\  0< {\rm Re} s < {1\over 4}.\eqno(3.5)$$
Substituting this value into (3.3), we end up with the estimate (3.1).  Moreover,  taking the Mellin transform (1.14) of both sides  in (1.2) under the condition $(G_\alpha g) (x) \in L_{\gamma,1}(\mathbb{R}_+)$, we employ (1.5) and change the order of integration via Fubini theorem in the right-hand side of the obtained equality.  Hence we obtain

$$\frac{\Gamma(1+\alpha-s)}{\Gamma(\alpha+s)\Gamma(1-s)} (G_\alpha  g)^* (s)$$

$$ =  {\Gamma(1/2-s) \over  \pi^{5/2}}  \int_{-\infty}^\infty  \Gamma(s+i\tau)\Gamma(s-i\tau) \cosh (\pi\tau)  g(\tau) d\tau. \eqno(3.6)$$
Hence inverting  the  Mellin transform  of both sides in (3.6) and using (2.7), we derive (3.2).    The convergence of the integral with respect to $s$ in (3.2) is absolute under the condition $|s|^{1/2-\gamma} e^{\pi|s|/2} (G_\alpha g)^*(s) \in L_1( \gamma - i\infty,\  \gamma +i\infty)$, and it can be easily verified,  recalling the Stirling asymptotic formula for the gamma-function.  Theorem 4 is proved. 

\end{proof} 

The inversion formula for the index transform (1.2) is given by

{\bf Theorem 5}.  {\it  Let  $\alpha \in \mathbb{R},\ g(z/i)$ be an even analytic function in the strip $D= \left\{ z \in \mathbb{C}:   |{\rm Re} z | \le \mu,  0< \mu < 1/4\right\} ,\  g(0)=0,  g^\prime(0)=0$ and $g(z/i) \in L_2\left(\mu-i\infty, \mu + i\infty\right)$.   Besides, if  $(G_\alpha g)^*(s)$ is analytic in $D$ and $|s|^{1/2-\nu} e^{\pi|s|/2} (G_\alpha g)^*(s) \in L_1( \nu - i\infty,\  \nu +i\infty)$ over any vertical line ${\rm Re} s=\nu$ in the strip,  then for all  $x \in \mathbb{R}$ the  following inversion formula holds for the index transform (1.2)} 
$$ g(x)  = -  \pi x \sinh(\pi x)  \left( \alpha^2 \int_0^\infty \int_{y}^\infty   {\rm Re} \left[ J_{\alpha +ix} \left(\sqrt t\right)  Y_{\alpha - ix} \left(\sqrt t\right) \right]  (G_\alpha g)(y) {dt dy\over t y}\right.$$

$$\left. +  \int_0^\infty {d\over dy} \left[   {\rm Re} \left[ J_{\alpha +ix} \left(\sqrt y\right)  Y_{\alpha - ix} \left(\sqrt y\right) \right]  \right] (G_\alpha g)(y) dy\right).\eqno(3.7)$$

\begin{proof}    Indeed,  since the integrand in the left-hand side of (3.2) is analytic in the strip $|{\rm Re} s| \le \mu$ and absolutely integrable there over any vertical line, we shift the contour to the left, integrating over $(\nu-i\infty, \nu+i\infty), -\mu < \nu < 0$.  Then multiplying both sides by $ e^{-y/2} K_{ix} \left({y/2} \right) y^{-1}$ and integrating  with respect to $y$ over $(0, \infty)$, we change the order of integration in the left-hand side of the obtained equality due to the absolute convergence of the  iterated integral.  Moreover,  appealing   to relation (8.4.23.3) in \cite{prud}, Vol. III, we calculate the inner integral  to find  the equality 
$$ {1\over 2\pi i}   \int_{\nu -i\infty}^{\nu  +i\infty} \frac {\Gamma(- s+ix)\Gamma(- s-ix) \Gamma(1+\alpha-s)}{\Gamma(1/2 -s) \Gamma(\alpha+s)\Gamma(1-s)} (G_\alpha g)^*(s) ds$$
$$ =  {1\over \pi^2\sqrt \pi } \int_0^\infty  K_{ix} \left({y\over 2} \right) \int_{-\infty}^\infty  K_{i\tau} \left({y\over 2} \right)  g(\tau) {d\tau dy\over y}. \eqno(3.8)$$
In the meantime the right-hand side of (3.8) can be written, employing  the  representation  of the Macdonald function in terms of the modified Bessel function of the first kind $I_z(y)$ \cite{erd}, Vol. II 

$${2 i \over \pi}  \sinh(\pi\tau)  K_{i\tau} \left({y\over 2} \right) =   I_{ -i\tau} \left({y\over 2} \right) -  I_{i\tau} \left({y\over 2} \right),$$
the substitution $z=i\tau$ and the property  $g(- z/i)=  g(z/i)$.  Hence it becomes

$${1\over \pi^2 \sqrt \pi } \int_0^\infty  K_{ix} \left({y\over 2} \right)  \int_{-\infty}^\infty  K_{i\tau} \left({y\over 2} \right) g(\tau) {d\tau dy\over y}$$

$$ = -  {1\over \pi i  \sqrt \pi } \int_0^\infty  K_{ix} \left({y\over 2} \right)  \int_{-i\infty}^{i\infty}   I_{ z} \left({y\over 2} \right)  g\left({z\over i}\right) {dz  dy\over y \sin(\pi z)}.\eqno(3.9)$$
On the other hand, according to our assumption $g(z/i)$ is analytic in $D$   and $|g(z/i) |^2$ is  integrable  over the vertical line $(\mu-i\infty, \mu + i\infty ), \mu >0$ in the strip.  Hence,  appealing to the inequality for the modified Bessel   function of the first  kind  (see \cite{yal}, p. 93)
 $$|I_z(y)| \le I_{  {\rm Re} z} (y) \  e^{\pi |{\rm Im} z|/2},\   0< {\rm Re} z \le  \mu,$$
one can move the contour to the right in the latter integral in (3.9) because the integrand is absolutely integrable. Then 

$$- {1\over \pi i \sqrt \pi } \int_0^\infty  K_{ix} \left({y\over 2} \right)  \int_{-i\infty}^{i\infty}   I_{ z} \left({y\over 2} \right)  g\left({z\over i}\right) {dz  dy\over y \sin(\pi z) }$$

$$=-  {1\over \pi i \sqrt \pi } \int_0^\infty  K_{ix} \left({y\over 2} \right)  \int_{\mu -i\infty}^{\mu+ i\infty}   I_{ z} \left({y\over 2} \right)  g\left({z\over i}\right) {dz  dy\over y\sin(\pi z) }.\eqno(3.10)$$
Hence one can  interchange the order of integration in the right-hand side of (3.10) due to the absolute and uniform convergence.  Then using  the value of the integral (see relation (2.16.28.3) in \cite{prud}, Vol. II)
$$\int_0^\infty K_{ix}(y) I_z(y) {dy\over y} = {1\over x^2 + z^2},$$ 
we find

$$ - {1\over \pi  i \sqrt \pi } \int_0^\infty  K_{ix} \left({y\over 2} \right)  \int_{\mu -i\infty}^{\mu+ i\infty}   I_{ z} \left({y\over 2} \right)  g\left({z\over i}\right) {dz  dy\over y \sin(\pi z) }$$

$$=  -   {1\over \pi i \sqrt \pi }  \int_{\mu -i\infty}^{\mu+ i\infty}  {g\left(z/ i\right) \over \left(x^2+ z^2\right) \sin(\pi z) } dz$$

$$ = - {1\over 2 \pi i \sqrt \pi }  \left( \int_{-\mu +i\infty}^{- \mu- i\infty}   +   \int_{\mu -i\infty}^{ \mu +  i\infty}   \right)  {  g(z/i) \  dz \over   z (z-ix) \sin(\pi z)}.\eqno(3.11)$$
Thus we are ready to apply the Cauchy formula in the right-hand side of the latter equality in (3.11) under conditions of the theorem.  Hence 

$$ - {1\over \pi  i \sqrt \pi } \int_0^\infty  K_{ix} \left({y\over 2} \right) \int_{\mu -i\infty}^{\mu+ i\infty}   I_{ z} \left({y\over 2} \right)  g\left({z\over i}\right) {dz  dy\over y \sin(\pi z) }$$

$$=  { 1\over \sqrt \pi}  \ { g(x) \over  x \sinh(\pi x)} ,\quad x  \in \mathbb{R} \backslash \{0\}.\eqno(3.12)$$
Combining with (3.8), we establish the equality 

$$   { 1\over \sqrt \pi}  \ { g(x) \over  x \sinh(\pi x)} =     {1\over 2\pi i} \int_{\nu -i\infty}^{\nu  +i\infty} \frac {\Gamma(- s+ ix)\Gamma(- s-ix) \Gamma(1+\alpha-s)}{\Gamma(1/2-s) \Gamma(\alpha+s)\Gamma(1-s)} $$

$$\times (G_\alpha g)^*(s) ds ,\quad x \in \mathbb{R} \backslash \{0\}.\eqno(3.13)$$
The right-hand side of (3.13) can be treated, in turn, appealing to the Parseval equality (1.15), conditions of the theorem and integral representation (2.26).   Therefore we derive

$$ g(x)  = -  \pi x\sinh(\pi x)  \int_0^\infty  {h(y) \over y}   \int_{y}^\infty   {\rm Re} \left[ J_{\alpha +ix} \left(\sqrt t\right)  Y_{\alpha - ix} \left(\sqrt t\right) \right]  {dt dy\over t},\ x \in \mathbb{R}, \eqno(3.14)$$
where (see (1.16))

$$h(y)=  {1\over 2\pi i} \int_{\nu -i\infty}^{\nu  +i\infty} (\alpha^2-s^2) (G_\alpha g)^*(s) y^{-s} ds$$

$$= \alpha^2  (G_\alpha g)(y) -  y \left({d\over dy} y {d\over dy}\right) (G_\alpha g)(y).\eqno(3.15)$$
Substituting the latter expression in (3.14), we integrate twice by parts, eliminating the integrated terms in order to arrive at the inversion formula (3.7). Theorem 5 is proved. 

 \end{proof}

\section{Boundary   value problem}

In this section we will employ  the index transform (1.2) to investigate  the  solvability  of a boundary  value  problem  for the following fourth    order partial differential  equation, involving the Laplacian

$$  r \left[ \left(  x {\partial  \over \partial x}  + y  {\partial   \over  \partial y} +2 \right)^2   - \alpha^2 \right]  \Delta u  $$
$$ +   \left(  x {\partial  \over \partial x}  + y  {\partial   \over  \partial y} \right)^2 u  +  {3\over 2}  \left(  x {\partial  \over \partial x}  + y  {\partial   \over  \partial y} \right) u   +  {u  \over 2 }  = 0, \  (x,y) \in  \mathbb{R}^2,\eqno(4.1)$$ 
where $r= \sqrt{x^2+y^2},\ \Delta = {\partial^2 \over \partial x^2} +  {\partial^2 \over \partial y ^2}$ is the Laplacian in $\mathbb{R}^2$.   In fact, writing  (4.1) in polar coordinates $(r,\theta)$, we end up with the equation  

$$r^2 {\partial^4 u \over \partial r^4} +    {\partial ^4  u \over \partial r^2  \partial \theta^2}+  6 r {\partial^3 u  \over \partial r^3} +  {1\over r} {\partial^3 u \over \partial r  \partial \theta^2}  +   \left( 7-\alpha^2 +r\right) {\partial^2  u \over \partial r^2}$$

$$ - {\alpha^2 \over r^2}   {\partial^2   u \over \partial \theta^2} + \left( {1-\alpha^2\over r}  + {5\over 2} \right) {\partial  u \over \partial r}  +  {u \over 2r }  = 0.\eqno(4.2)$$

{\bf Lemma 3.} {\it  Let $g(\tau)  \in L_2\left(\mathbb{R};  \tau^2  d\tau\right),\  \beta \in (0, 2\pi)$. Then  the function
$$u(r,\theta)=    \int_{-\infty}^\infty   {\mathcal W}_\alpha(r,\tau)\  {\sinh(\theta \tau)\over \sinh(\beta\tau) }   g(\tau) d\tau,\eqno(4.3)$$
 satisfies   the partial  differential  equation $(4.2)$ on the wedge  $(r,\theta): r   >0, \  0\le \theta \le  \beta$, vanishing at infinity.}

\begin{proof} The proof  is straightforward by  substitution (4.3) into (4.2) and the use of  (1.11).  The necessary  differentiation  with respect to $r$ and $\theta$ under the integral sign is allowed via the absolute and uniform convergence, which can be verified  analogously to estimates (3.3) under the  condition $g \in L_2\left(\mathbb{R};  \tau^2 d\tau\right)$.   Finally,  the condition $ u(r,\theta) \to 0,\ r \to \infty$  is due to inequality (2.5). 
\end{proof}

Finally  we will formulate the boundary  value problem for equation (4.2) and give its solution.

{\bf Theorem 6.} {\it Let  $g(x)$ be given by formula $(3.7)$ and its transform $(G_\alpha g) (t)\equiv G_\alpha (t)$ satisfies conditions of Theorem 5.  Then  $u (r,\theta),\   r >0,  \  0\le \theta \le \beta$ by formula $(4.3)$  will be a solution  of the boundary  value problem for the partial differential  equation $(4.2)$ subject to the conditions}
$$u(r,0) = 0,\quad\quad   u(r,\beta) =  G_\alpha (r).$$

\bigskip
\centerline{{\bf Acknowledgments}}
\bigskip

The work was partially supported by CMUP (UID/MAT/00144/2013),  which is funded by FCT (Portugal) with national (MEC),   European structural funds through the programs FEDER  under the partnership agreement PT2020, and Project STRIDE - NORTE-01-0145-FEDER- 000033, funded by ERDF - NORTE 2020. 

\bibliographystyle{amsplain}

\end{document}